\documentstyle[12pt]{article}

\newcommand{\be}{\begin{equation}}
\newcommand{\ee}{\end{equation}}
\newcommand{\bn}{\begin{eqnarray}}
\newcommand{\en}{\end{eqnarray}}
\newcommand{\bd}{b^\dagger}
\newcommand{\ad}{a^\dagger}

\normalbaselines
\begin{document}
\begin{titlepage}
\begin{flushright}
{\bf HU-SEFT \ R \ 1995-09}
\end{flushright}
\begin{flushright}
{\bf CBPF-NF-022/95}
\end{flushright}
\begin{center}

\vspace*{1.0cm}

{\Large{\bf Path Integrals with Generalized Grassmann Variables}}

\vskip 1.5cm

{\large {\bf M. Chaichian}}

\vskip 0.5cm

High Energy Physics Laboratory, Department of Physics \\
and Research Institute for High Energy Physics, \\
University of Helsinki\\
P.O.Box 9 (Siltavuorenpenger 20 C), FIN-00014,  \\
Helsinki, Finland

\vskip 0.2cm

$and$

\vskip 0.2cm

{\large {\bf A.P.Demichev}}\renewcommand{\thefootnote}
{\dagger}\footnote{on leave of absence from
Nuclear Physics Institute,
Moscow State University,
119899, Moscow, Russia}

\vskip 0.5cm

Centro Brasileiro de Pesquisas Fisicas - CBPF/CNPq, \\
Rua Dr.Xavier Sigaund, 150, \\
22290-180, Rio de Janeiro, \\
RJ Brasil

\end{center}

\vspace{3 cm}

\begin{abstract}
\normalsize

We construct path integral representations for the evolution
operator of q-oscillators with root of unity values of q-parameter using
Bargmann-Fock representations
with commuting and non-commuting variables, the differential calculi being
q-deformed in both cases. For $q^2=-1$
we obtain a new form of Grassmann-like path integral.

\end{abstract}
\end{titlepage}

\section{Introduction}

Path integrals over Grassmann variables play important role in modern field
theory. In the simplest case of one degree of freedom it describes the quantum
mechanical system with two possible states. It is natural to generalize the
consideration for the case of quantum mechanical system with arbitrary but
finite number of states. As is well known, operator description of such systems
is realized with help of q-oscillators algebra \cite{qOsc}, parameter q being
root of unity: $q^p=1$ for some integer $p$. Note that in general  quantum
algebras in the root of unity case exhibit rich representation behaviours  with
very special properties (including bound on their dimensionality) and are
important in two-dimensional conformal field theories \cite{A-G} and in
statistical mechanics models \cite{Vega}. We shall consider even integer $p$,
so that $q^{(k+1)}=-1$ for $k=p/2-1$. The case $q^2=-1$ corresponds to the
usual fermionic oscillator. Of      course, generalization to field theory
requires the consideration of more complicated algebra of system of
q-oscillators and its precise physical meaning depends on the chosen form of
commutation relations for the different q-oscillators (bosonic excitations on
lattice with finite number of states at a given site or system with fractional
statistics). Here we consider one degree of        freedom only and the aim of
this letter is to present possible
Bargmann-Fock representations of the q-oscillators with finite Fock space and
to    construct corresponding q-deformed path integral. An interesting attempt
to acheive this was made in \cite{BF} but as we discussed earlier \cite{CD} the
authors did not take into account that the creation and destruction operators
used by them are not conjugated to each other in the root of unity case and as
the final result they obtained discrete approximation for the path integral
representation of (non-unitary) evolution operator.

If the Hilbert space of a system corresponds to representation of some Lie
group the corresponding path integral can be constructed with help of
generalized coherent states (see, e.g., \cite{IKG} and refs. therein). For our
construction it is not necessary to assume that there is any symmetry. In
general we shall follow the method suggested by us in \cite{CD} for the
q-deformed path integral in the case of real q-parameter. It is interesting
that even in the simplest Grassmann-like case $q^2=-1$ there are different
forms of BF representation and path integrals depends on the way how the Planck
constant $\hbar$ enters the commutation relation. To generalize the usual
Grassmann path integral for the case of higher roots we suggest another algebra
of variables on different time slices, the Grassmann integral being independent
of the modification.

Such path integrals can be useful for consideration of systems on a lattice
\cite{Kul}, in anyon physics \cite{Ler}, for description of quantum particles
on finite sets of points \cite{DM-H} and for q-deformed string theory
\cite{CGG}.

\section{Path integrals for different choices of variables}

The usual starting commutation relations (CR) for
q-oscillator operators is \cite{qOsc}
$$
 a\ad - q \ad a = q^{- N} \ , \qquad
\ad a = [N]_q\equiv \frac{q^N-q^{-N}}{q-q^{-1}}\ ,     $$
where $q^{k+1}=-1$ for some integer $k$ and $\ad$ is hermitian conjugated to
$a$. To construct the Bargmann-Fock (BF) representation we introduce the
operators
$b$ and $\bar b$ (well known in the case of real q-parameter)
$$
b=q^{N/2}a\ ,\qquad \bar b = \ad q^{N/2}           $$
with CR
\be
b\bar b -q^2\bar b b=1\ ,                          \label{1}
\ee
\be
b\bar b -\bar b b = q^{2N} \ ,                     \label{2}
\ee
$$
\bar b b=[N]\equiv\frac{q^{2N}-1}{q^2-1}  \ .
$$
For root of unity value of parameter $q$ operator $\bar b$ is not hermitian
conjugate to $b$, but instead
$$
\bd=\bar b q^{-N}, \qquad   \bar \bd=q^{-N}b \ .     $$
Fock space representation of the operators is the following
\be
b|n\rangle=q^{n-1}\sqrt{[n]_q}|n-1\rangle \ , \quad
\bar b|n\rangle=\sqrt{[n+1]_q}|n+1\rangle \ , \quad
N|n\rangle=n|n\rangle \ .                           \label{3}
\ee
As $[k+1]_q=0$, this Fock space is k-dimensional.

There are two possibilities for explicit realization of (\ref{3}). One can
construct BF representation in the space of antiholomorphic functions with the
basis
$$
\psi_n=\frac{\bar z^n}{\sqrt{[n]_q!}}\ ,                   $$
with either commuting $\bar z,z$ variables or non-commuting ones.
In both cases the variables must satisfy the condition of nilpotence
$
z^{k+1}=\bar z^{k+1}=0                           $
to provide the finiteness of the Fock space.

We shall consider this  on the example of Grassmann-like $q^2=-1$ case
\renewcommand{\thefootnote}{1}\footnote {We call this case Grassmann-{\it like}
because though operators $b$ and $\bar b$ satisfy the anticommutation
relations, we use for the construction of the path integrals
either commuting variables or variables with non-Grassmann algebra on
different time slices.}.
Physicaly the possibilities correspond to the different forms of the
CR (\ref{1}),(\ref{2}) after recovering of Planck constant $\hbar$.
For $q^2=-1$ case and in representation (\ref{3}) both CR (\ref{1}) and
(\ref{2}) become the same
\be
b\bar b +\bar b b=1\ .                                 \label{4}
\ee
As usual one understands that Planck constant in this CR is equal to unity. But
it is possible to recover it in two different ways. One way is to present
(\ref{4}) in the form
\be
[b,\bar b]=\hbar (1-2\gamma\bar b b)\ ,                 \label{5}
\ee
where $\gamma$ is dimensionful constant such that $\hbar\gamma=1$. This CR
corresponds to a curved phase space dynamics with Poisson braket
\be
\{z,\bar z\}_P=i(1-2\gamma \bar z z)\equiv i\omega^{-1} (\bar z z) \ ,
                                              \label{6}
\ee
where $z,\bar z$ are the classical counterparts of $b,\bar b$ and
$\omega=\omega(\bar z z)dz\wedge d\bar z$ is the
symplectic form. Evolution operator corresponding to such quantization
must be expressed with help of q-path integral with commuting variables.

The key observation for the developing of path integral representations is that
CR written in the form (\ref{1}) does not lead
unavoidably to BF variables with CR
$
z\bar z=q^2\bar z z\ ,                   $
but to CR
\be
\bar \partial\bar z - q^2 \bar z\bar \partial=1  \label{6a}            \ee
only. So one can consider two parametric differential calculus
on the quantum plane \cite{KR} and for the case of curved phase space
(\ref{5}) choose the commuting complex variables $z\bar z=\bar z z$ with the
following nontrivial CR only
\be
\bar\partial\bar z+\bar z\bar\partial=1\ , \qquad
d\bar z\bar\partial=-\bar\partial d\bar z\ ,  \qquad
\bar zd\bar z = -(d\bar z) \bar z                  \label{7}
\ee
with their conjugated counterparts.
The  creation $\bar b=\bar z$ and destruction
$b=\bar\partial$ operators act in the BF space with the basis $\{\psi_0=1,
\psi_1=\bar z\}$ and the scalar product
$$
\int d\bar z dz e^{\bar z z}\bar\psi_n\psi_m = \delta_{nm}  \ ,    $$
where the integral is defined by the usual Berezin rules
$$
\int d\bar z\ \bar z=\int dz\ z= 1\ ,\qquad
\int d\bar z=\int dz =0\ .                                        $$
The method of the path integral construction is analogous to
that in the case of usual Grassmann path integral (see e.g. \cite{Ber}).
As a result the evolution operator kernel takes the form
$$
U(t^{\prime\prime}-t^\prime)= \int\left(\prod_t
\frac{d\bar z(t)dz(t)}{1-2\bar z(t)z(t)}\right)
\exp\left\{\bar z(t^{\prime\prime}) z(t^{\prime\prime})
-\int^{t^{\prime\prime}}_{t^\prime}\left(
\bar z(t)\dot z(t)+iH(\bar z(t)z(t))\right) dt\right\}\ ,      $$
with $\dot z(t)$ denoting the usual derivative with respect to time:
$\dot z(t)=dz(t)/dt$.
Note that the integral measure corresponds to the form of nontrivial
Poisson bracket (\ref{6}) after putting $\gamma=\hbar=1$.

If after the recovering of Planck constant the CR (\ref{4}) takes the form
$$
b\bar b +\bar b b=\hbar\ ,
$$
then it corresponds to Grassmann phase space with anticommuting variables.
 In this
case the BF variables seems to be identical to the usual Grassmann path
 integral.
But there is one subtlety. In the usual construction not only $z(t_i)$
and $\bar z(t_j)$ anticommute for any time slices $t_i,t_j$ but the
same variables on different time slices anticommute also, e.g.
$z(t_i)z(t_j)+z(t_j)z(t_i)=0$. Such CR can not be generalized to other
roots of unity values of parameter $q$. So we introduce another CR
on different time slices (cf. \cite{CD})
\be
\begin{array}{ccc}
z(t_i)\bar z(t_j) + \bar z(t_j) z(t_i)=0\ ,& \quad &
z^2(t_i)=\bar z^2(t_i) = 0 \ , \\
\bar z(t_i)\bar z(t_j)=\bar z(t_j)\bar z(t_i)\ ,& \quad &
z(t_i)z(t_j)=z(t_j)z(t_i)\ .
\end{array}                                                \label{11}
\ee
One can check that for such CR all ingredients of path integral
construction (scalar product measure, relation between normal symbol
and kernel of operators, etc.) remain the same as in the case of
usual Grassmann path integral. So the result
proved to be the same also.

The CR (\ref{11}) can be easily generalized to other values of deformation
parameter, for example, for $q^3=-1$ one has
\be
\begin{array}{ccc}
z(t_i)\bar z(t_j)=q^2\bar z(t_j)z(t_i)\ ,& \quad &
z^3(t_i)=\bar z^3(t_i) = 0 \ , \\
\bar z(t_i)\bar z(t_j)=\bar z(t_j)\bar z(t_i)\ ,&\quad &
z(t_i)z(t_j)=z(t_j)z(t_i)\ .
\end{array}                                                \label{12}
\ee
In this case the BF representation is defined on the space of function
with the basis

\be
\psi_0=1\ ,\quad \psi_1=\bar z\ , \quad \psi_2=\bar z^2\ ,\label{12a}
\ee
which is orthonormal with respect to the scalar product

\be
\int d\bar z dz\  \overline{\psi_n(\bar z)}
\mu(\bar z z)\psi_m(\bar z) = \delta_{nm}\ ,    \label{12b}
\ee
where
$$
\mu(\bar z z)=1+q^2 \bar z z+q^2 \bar z^2 z^2 \ .           $$
Here the integral is defined by the relation \cite{BF}:
\be
\int d\bar z dz\  z^n\bar z^m = \delta_{nk}\delta_{mk}[k]_q! \ ,
\quad n,m=0,1,2\ .
                                                             \label{12c}
\ee

Because of the nilpotence a general Hamiltonian has the form
\be
H=\omega (u\bar bb+v\bar b^2b^2)   \ ,                      \label{12d}
\ee
where the constants $u,v$ are restricted by the hermiticity condition
$H^\dagger=H$ and can take three couples of the values:
{\it i} )$u=1,\ v=-q$; {\it ii})$u=1,\ v=1-2q$; {\it iii})$u=0,\ v=-q^2$.
In this case and in the cases of higher roots of unity the measure of scalar
product in BF spaces and relation between normal symbols and kernels of
operators are not expressed with help of appropriate q-deformed exponents. This
prevents from writing the general expression of path integral for arbitrary
root of unity. But due to nilpotence for any given value of parameter $q$ it
can be done analogously to the case of real q-parameter \cite{CD}.

As usual the action of any operator $A$ in BF Hilbert space can be
represented with the help of its kernel ${\cal A}$
\begin{equation}
(Af)(\bar z_1) = \int d\bar z_2 dz_2\ {\cal A}(\bar
z_1,z_2)f(\bar z_2)\ ,                            \label{26}
\end{equation}
where
\begin{equation}
{\cal A}(\bar z_1,z_2)=\sum_{m,n=0}^2 A_{mn}
\bar z_1^m z_2^n\ .                               \label{27}
\end{equation}
Here one more pair of q-commuting coordinates is introduced, the CR for
different pairs being defined by (\ref{12}). Now we can express $A_{mn}$
through scalar product
\begin{equation}
A_{mn}=<\psi _m \| A \|\psi _n>\                  \label{29}
\end{equation}
and find a kernel of any operator by direct calculation.
We consider the usual
Schr{\"o}dinger equation (cf. discussion in \cite{CD})
\begin{equation}
i\frac{d\ }{dt} \Psi (\bar z,t)=H(\bar b,b)\Psi(\bar z,t)\ ,
                                                  \label{36}
\end{equation}
with the Hamiltonian (\ref{12d}).
The integral kernel for the infinitesimal operator
$$
 U\approx 1-iH\Delta t\ , $$
takes the form
\begin{equation}
U(\bar z, z;\Delta t) \approx g(\bar z z)exp\{-iH_{eff}\Delta t\}
\ ,                                                  \label{38}
\end{equation}
where $g(\bar z z)$ is the kernel of identity operator
$$
g(\bar z z)=\sum_n\psi_n(\bar z)\overline{\psi_n(\bar z)} = 1+\bar z z +
\bar z^2 z^2 \ ,
$$
and  effective Hamiltonian $H_{eff}$ is defined by the relation
\be
H_{eff}(\bar z z)=g^{-1}(\bar z z){\cal H}(\bar z z)
\ ,                                             \label{13a}
\ee
where ${\cal H}(\bar z z)$ is kernel of the initial Hamiltonian.
 Using the nilpotence one obtains from (\ref{13a})
$$
H_{eff}(\bar z z)=\omega (u\bar zz+q(u+v-qu)\bar z^2z^2)   \ .   $$

Now we can write the convolution of $K$ infinitesimal evolution
operator kernels
$$
U(\bar z_Kz_{K-1})*U(\bar z_{K-1}z_{K-2})*...*U(\bar z_1z_0)\  =  $$
$$
=\int d\bar z_{K-1}dz_{K-1}...d\bar z_1dz_1\ \mu (\bar z_{K-1}z_{K-1})
...\mu (\bar z_1z_1)                         $$
\begin{equation}
\times g(\bar z_Kz_{K-1})...
g(\bar z_1z_0)e^{-iH_{eff}(\bar z_Kz_{K-1})\Delta t}
...e^{-iH_{eff}(\bar z_1z_0)\Delta t}\ .             \label{40}
\end{equation}
Due to nilpotence it is possible to convert the functions $g$ and $\mu$
to the exponents and in the continuous limit
$\Delta t\rightarrow 0$ one obtains the path integral
for the case under consideration ($q^3=-1$)
\begin{eqnarray}
U(t^{\prime\prime}-t^\prime)= \int\left(\prod_t
d\bar z(t)dz(t)(1+q\bar z(t)z(t))\right)
(1+\bar z(t^{\prime\prime}) z(t^{\prime\prime})+
\bar z^2(t^{\prime\prime}) z^2(t^{\prime\prime}))  \nonumber \\
\times \exp\left\{-\int^{t^{\prime\prime}}_{t^\prime}\left[
(1+(1+2q)\bar z z)\bar z(t)\dot z(t)+
iH_{eff}(\bar z(t)z(t))\right] dt\right\}
\ .                                                    \label{13}
\end{eqnarray}

The path integrals for higher roots can be derived in the same way. For
example, for the next value $q^4=-1$ one has
\begin{eqnarray}
U(t^{\prime\prime}-t^\prime)= \int\left(\prod_t
d\bar z(t)dz(t)\left[1+(1+i)\bar z(t) z(t) +(2c-1)\bar z^2(t) z^2(t) +
\right.\right. \nonumber\\
\left.\left. i(1-c)\bar z^3(t)z^3(t)\right]\right)\exp\{\bar
z(t^{\prime\prime}) z(t^{\prime\prime})+
(c-i/2)\bar z^2(t^{\prime\prime}) z^2(t^{\prime\prime})+(2c-i/3)\bar
z^3(t^{\prime\prime}) z^3(t^{\prime\prime})\} \nonumber  \\
\times \exp\left\{-\int^{t^{\prime\prime}}_{t^\prime}\left[
(1-(1/2+ic)\bar z z-3(2c-i/2)\bar z^2z^2)\bar z(t)\dot z(t)+
iH_{eff}(\bar z(t)z(t))\right] dt\right\}\ ,    \label{13aa}
\end{eqnarray}
where $c\equiv 2^{-1/4}=1/\sqrt{[2]_q}$ and
$H_{eff}$ is defined again by the relation (\ref{13a}) with the corresponding
kernel
$$
g(\bar z z)=1+\bar z z +c\bar z^2 z^2+c\bar z^3z^3\ .
$$

The kernels $U$ in Eqs. (\ref{13}) and (\ref{13aa}) correspond to the product
of
unitary up to a
$(\Delta t)^2$-terms operators with subsequent limit $\Delta t
\rightarrow 0$. So these kernels seem to be the kernels of unitary evolution
operators. However, the integrands of these path integrals have not the form
$\exp\{iS\}$, where $S$ is a real functional. This is quite unusual and
prevents
from ordinary interpretation of a path integral. The obvious reason for this is
the CR (\ref{12}) for $z$ and $\bar z$ coordinates which contain complex
parameter $q$.

Let us consider for $q^3=-1$ another possibility using commuting variables with
CR (\ref{6a}) for variables and derivatives in analogy with Grassmann-like
case. The basis of BF representation has the same form (\ref{12a}) as in the
case of noncommuting variables and is orthonormal with respect to scalar
product (\ref{12b}) with the mesure $\mu=1+\bar z z+\bar z^2 z^2$, the integral
being defined by (\ref{12c}) (but now with commuting variables). The derivation
of the path integral is essentially the same as in the case of noncommuting
variables and the result is
\begin{eqnarray}
U(t^{\prime\prime}-t^\prime)= \int\left(\prod_t
d\bar z(t)dz(t)(1+2\bar z(t)z(t)+3\bar z^2(t) z^2(t))\right) (1+\bar
z(t^{\prime\prime}) z(t^{\prime\prime}) \nonumber \\
+\bar z^2(t^{\prime\prime}) z^2(t^{\prime\prime}))
\exp\left\{-\int^{t^{\prime\prime}}_{t^\prime}\left[
(1+\bar z z)\bar z(t)\dot z(t)+
iH_{eff}(\bar z(t)z(t))\right] dt\right\} \ . \label{13aaa}
\end{eqnarray}
The effective Hamiltonian is defined by the general formula (\ref{13a}) again,
but now for all hermitian Hamiltonians of the form (\ref{12d}) listed above
$H_{eff}$ proved to be real function: {\it i}) $H_{eff}=\omega\bar z z$ for
$u=1,v=-q$; {\it ii}) $H_{eff}=\omega (\bar z z-\bar z^2z^2)$ for $u=1,v=1-2q$;
{\it iii}) $H_{eff}=\omega\bar z^2 z^2$ for $u=0,v=-q^2$. Note that the normal
symbols and the kernels of all hermitian Hamiltonians are complex. Thus in the
case of commuting variables $\bar z$ and $z$ the path integral has the usual
form  $\exp \{iS\}$ of integrand with real functional $S$ even for complex
parameter $q$ of an oscillator deformation.

After  restoring the Planck constant $\hbar$ the
expression (\ref{13aaa}) takes the form (for definiteness we consider the
Hamiltonian (\ref{13a}) with $u=1, v=-q$)
\begin{eqnarray}
U(t^{\prime\prime}-t^\prime)= \int\left[\prod_t
\frac{d\bar z(t)dz(t)}{\hbar}\left( 1+2\frac{\bar z(t)z(t)}{\hbar}+3\frac{\bar
z^2(t) z^2(t)}{\hbar^2}\right)\right] \left( 1+\frac{\bar z(t^{\prime\prime})
z(t^{\prime\prime})}{\hbar} \right.
\nonumber \\
\left. +\frac{\bar z^2(t^{\prime\prime}) z^2(t^{\prime\prime})}{\hbar^2}\right)
\exp\left\{-\frac{1}{\hbar}
\int^{t^{\prime\prime}}_{t^\prime}\left[
\left( 1+\frac{\bar z z}{\hbar}\right)\bar z(t)\dot z(t)+
i\omega\bar z(t)z(t)\right] dt\right\} \ . \label{13aaaa}
\end{eqnarray}
One can see that  there is simply no quasiclassical approximation
corresponding to the limit $\hbar\rightarrow0$
(since the terms in the measure and in the action diverge). This is an expected
result because, as it is well known, the quasiclassical approximation for
spin-like systems considered here corresponds to the limit
$S\rightarrow\infty$, where S is the spin of the system \cite{Ber2}. Systems
with
fixed spins (numbers of states) have not (quasi)classical limit, as is
confirmed once again by the expression (\ref{13aaaa}).

  Let us notice that while the real phase space variables (coordintaes and
momenta)
are not suitable in describing spin systems, the BF basis is. The simplest and
well known example is the Grassmann algebra variables for describing spin 1/2
particles. It is also of
interest to mention that in a recent attempt to formulate a regularized
quantum field theory by requiring nonzero minimal uncertainties in the
positions and in the momenta, for the real q case the use of Bargmann-Fock
space
appears to be the appropriate one \cite{Kempf}. The q-deformed path integral
 previously given
\cite{CD} for the real q values, has also been
formulated in the BF basis emerged in a natural way.

The fact that creation and destruction operators $\bar b,\ b$ are not hermitian
conjugate is quite unusual and technically inconvenient (but as they are
non-conjugate up to pure phase only this does not cause real problems). One can
try to use hermitian conjugated operators $\ad ,a$ for the path integral
construction noticing that the relation (\ref{1}) is in fact a polynomial one
due to nilpotence. In particular, in $q^3=-1$ case and in the Fock space
representation the relation (\ref{1}) is equivalent to the CR
\be
a\ad+(\ad )^2 a^2=1\ .                                    \label{15}
\ee
This CR looks like a natural generalization of the usual relation (\ref{4}) for
fermionic operators\renewcommand{\thefootnote}{2}\footnote {One can show that
the algebra of operators $\alpha\alpha^\dagger+(\alpha^\dagger)^n\alpha^n =1$,
has $(n+1)$-dimensional Fock representation; a general consideration of
polynomial algebras with finite-dimensional Fock representations will be given
elsewhere.}. To construct the BF representation for (\ref{15}) one must
 introduce
unusual "differential" operators ${\cal D},\bar{\cal D}$ with the properties
\be
\begin{array}{ccc}
\bar{\cal D}\bar z+\bar z^2\bar{\cal D}^2=1\ , & \quad & {\cal D} z+
z^2{\cal D}^2=1\ ,\\
\bar z^3=z^3=1\ . & &
\end{array}                                                   \label{16}
\ee
Consider for  definiteness the commuting variables $z\bar z=\bar z z$ (one can
consider other possibilities, e.g. anticommuting $z$ and $\bar z$). The basis
of the functions (\ref{12a}) is
orthonormal with respect to the scalar product (\ref{12b}) with function
$\mu = 1+\bar z z+\bar z^2 z^2 $, the integral being again defined by
(\ref{12c}). The operators $\bar z$ and $\bar{\cal D}$ are hermitian conjugate
with respect to this scalar product. One can check that the
normal symbols ${\cal A}_N (\bar z z)$ of operators $A$
and their kernels ${\cal A}(\bar z z)$ are related with help of the same
function
$\mu(\bar z z)$ that defines the scalar product, the kernel being defined in
the usual way
$$
A\psi(\bar z)=\int d\bar z^\prime dz^\prime {\cal A}(\bar z,z^\prime )
\psi (\bar z^\prime)\ .                                           $$
Thus using the nilpotence and usual procedure for path integral derivation one
obtaines
\begin{eqnarray*}
U(t^{\prime\prime}-t^\prime)= \int\left(\prod_t
d\bar z(t)dz(t)(1+2\bar z(t)z(t)+3\bar z^2(t)z^2(t))\right)
\mu(\bar z(t^{\prime\prime}) z(t^{\prime\prime})) \\
\times \exp\left\{-\int^{t^{\prime\prime}}_{t^\prime}\left[
(1+\bar z z)\bar z(t)\dot z(t)+
iH_{eff}(\bar z(t)z(t))\right] dt\right\}
\ .
\end{eqnarray*}
Analogous expressions can be derived for the algebras which correspond to Fock
spaces of arbitrary dimension.

\section{Conclusion}

We have succeeded in formulating q-deformed path integrals for root of unity
values of q-parameter for simplest
nontrivial quantum mechanical systems. Using different differential
and integral calculi on quantum plane we constructed the q-deformed
Bargmann-Fock representations and path integrals for commuting and
non-commuting variables. As we have shown different forms of path integrals are
possible even for the well known case of anticommuting creation and destruction
operators leading to different Grassmann-like path integrals. The existence of
non-standard forms of path integral for fermionic operators permited us to
generalize it to the case of q-oscillators with higher root of unity
q-parameter. The
essential distinction of the generalized Grassmann path integrals is that they
have a non-trivial integral measure and a non-Gaussian integrand.

\vspace{0.3cm}

{\bf Acknowledgements}

A.D.'s work was partially supported by the INTAS-93-1630 grant.

\end{document}